\documentclass[10.5pt]{amsart}
\usepackage{hyperref}
\usepackage{amscd,amsfonts,amssymb,amsmath}
\usepackage{hyperref}
\usepackage{url}
\usepackage{epsfig}
\usepackage{mathtools}
\usepackage{tabu}
\usepackage{tikz-cd}
\usepackage{doi}
\usepackage{siunitx}
\newtheorem{theoreme}{Theorem}[section]
\newtheorem{corollaire}[theoreme]{Corollary}
\newtheorem{definition}[theoreme]{Definition}
\newtheorem{lemme}[theoreme]{Lemma}
\newtheorem{proposition}[theoreme]{Proposition}

\newtheorem{remark}[theoreme]{Remark}
\numberwithin{equation}{subsection}
\newtheorem{thmA}{Theorem}

\newtheorem*{ack}{Acknowledgements}
\usepackage[all,cmtip]{xy}
\usepackage{graphicx}

\newcommand{\Id}{\operatorname{Id}}

\newcommand{\R}{\operatorname{R}}

\newcommand{\Z}{\operatorname{Z}}
\newcommand{\N}{\operatorname{N}}

\newcommand{\myrightleftarrows}[1]{\mathrel{\substack{\xrightarrow{#1} \\[-.9ex] \xleftarrow{#1}}}}
\begin{document} 
\title{Combinatorics of affine cactus groups}

\author{Hugo Chemin}
\address{Normandie Univ., UNICAEN, CNRS, LMNO, 14000 Caen, France}

\maketitle

\begin{abstract}
This article deals with the study of affine cactus groups from a combinatorial point of view. Those groups are extensions of cactus groups, which are related to braid and diagram groups and have gained an important place in many mathematics topics. We first show that affine cactus groups may be described as cactus groups on Coxeter groups of type $\widetilde{A}_n$. Then, we prove that these groups embed into a semi-direct product of Coxeter groups, which allows us to obtain a number of combinatorial properties of affine cactus groups, such as the solubility of the world problem or the fact that their centre is trivial.
\end{abstract}

\section{Introduction}
Affine cactus groups $AJ_n$ first appeared in the work of Ilin-Kamnitzer-Li-Przytycki-Rybnikov \cite{ilin_moduli_2023} as an extension of cactus groups. They show that there is an action of $\Z/n\Z$ on $AJ_n$ and that the group $\widetilde{AJ}_n=AJ_n \rtimes \Z/n\Z$ is the $S_n$-equivariant fundamental group of $\overline{M}_{n+2}^\sigma(\R)$, a real twisted form of $\overline{M}_{n+2}(\R)$, which is the Deligne-Mumford compactification of real stable curves of genus $0$ with $n+2$ marked points \cite{ilin_moduli_2023}.\\
Affine cactus groups are also trickle groups, which are groups generalizing right-angled Artin and Coxeter groups, as well as cactus groups \cite{trickle_2024}.\\
In the same way as for affine braid groups whose elements are braids on a circle (\cite{bellingeri_bodin_2016}), the elements of affine cactus groups are cacti on a circle (Fig \ref{figure1}). To express this circularity, Ilin-Kamnitzer-Li-Przytycki-Rybnikov \cite{ilin_moduli_2023} introduce the notion of circular intervals as follows : for $n \geq 2$ and $1 \leq i \neq j \leq n$, we define $[i,j]_c$ to be the set $\{i <_c i+1 <_c \ldots <_c j \}$ where $<_c$ is the cyclic order on $\Z/n\Z$. We say that $[i,j]_c$ is contained in $[k,l]_c$, or that $[i,j]_c$ is a sub-interval of $[k,l]_c$, denoted by $[i,j]_c \subset_c [k,l]_c$, if the elements of $[i,j]_c$ belong to $[k,l]_c$ with preservation of the order. For example, the circular interval $\{1 <_c 2 <_c 3 \}$ is included in $\{1 <_c 2 <_c 3 <_c 4 \}$ but not in $\{  3 <_c 4 <_c 1 <_c 2 \}$.\\
Let $n \geq 2$. Formally, the affine cactus group $AJ_n$ is generated by $\{\sigma_{i,j}, \ \ 1 \leq i \neq j \leq n\}$ subject to the following relations:
\begin{eqnarray}\label{Canonical-Presentation}
\sigma_{i,j}^{2} &=&1 \hspace*{5mm} \textrm{for } 1 \leq i \neq j \leq n,\\ 
\sigma_{i,j}\sigma_{k,l} &=& \sigma_{k,l}\sigma_{i,j} \hspace*{5mm} \textrm{if } [i,j]_c \cap [k,l]_c= \emptyset,\\
\sigma_{i,j}\sigma_{k,l} &=& \sigma_{k,l}\sigma_{s_{k,l}(j),s_{k,l}(i)} \hspace*{5mm} \textrm{if } [i,j]_c \subset_c [k,l]_c.
\end{eqnarray}
\definecolor{ttqqqq}{rgb}{0.2,0.,0.}
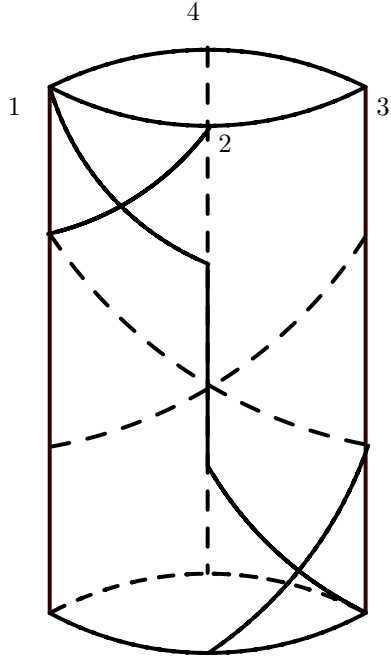
\begin{figure}[h]
\begin{center}
\begin{tikzpicture}[line cap=round,line join=round,x=0.7cm,y=0.7cm,rounded corners=3pt,line width=1.5pt]
\draw [color=ttqqqq] (2.,10.)-- (2.,0.);
\draw [color=ttqqqq] (8.,0.)-- (8.,10.);
\draw [shift={(5.,5.54076)}]  plot[domain=4.216139969888791:5.208637990880589,variable=\t]({1.*6.3007952972303425*cos(\t r)+0.*6.3007952972303425*sin(\t r)},{0.*6.3007952972303425*cos(\t r)+1.*6.3007952972303425*sin(\t r)});
\draw [shift={(7.0039759689188195,11.448961289710656)}]  plot[domain=3.423446006517393:4.320487491849233,variable=\t]({1.*5.209535902227472*cos(\t r)+0.*5.209535902227472*sin(\t r)},{0.*5.209535902227472*cos(\t r)+1.*5.209535902227472*sin(\t r)});
\draw [shift={(0.7972903234867222,12.327265862177464)}]  plot[domain=4.937100148799332:5.652155560798254,variable=\t]({1.*5.25557637836749*cos(\t r)+0.*5.25557637836749*sin(\t r)},{0.*5.25557637836749*cos(\t r)+1.*5.25557637836749*sin(\t r)});
\draw [dash pattern=on 7pt off 7pt] (5.,10.75)-- (5.,0.76292);
\draw [shift={(0.47524531358222666,11.946667214108514)},dash pattern=on 7pt off 7pt]  plot[domain=4.88725132960931:5.718004983249195,variable=\t]({1.*8.919056921582913*cos(\t r)+0.*8.919056921582913*sin(\t r)},{0.*8.919056921582913*cos(\t r)+1.*8.919056921582913*sin(\t r)});
\draw [shift={(9.755996962648146,12.327265862177464)},dash pattern=on 7pt off 7pt]  plot[domain=3.72470975717383:4.526076710533615,variable=\t]({1.*9.291397170904004*cos(\t r)+0.*9.291397170904004*sin(\t r)},{0.*9.291397170904004*cos(\t r)+1.*9.291397170904004*sin(\t r)});
\draw [shift={(5.,-5.61942)},dash pattern=on 7pt off 7pt]  plot[domain=1.0804267790841215:2.061165874505672,variable=\t]({1.*6.370077011810768*cos(\t r)+0.*6.370077011810768*sin(\t r)},{0.*6.370077011810768*cos(\t r)+1.*6.370077011810768*sin(\t r)});
\draw (5.,6.6)-- (5.,2.8);
\draw [shift={(11.0149001831839,6.266964312156488)}]  plot[domain=3.6644754198509464:4.2639930381990245,variable=\t]({1.*6.942540295556983*cos(\t r)+0.*6.942540295556983*sin(\t r)},{0.*6.942540295556983*cos(\t r)+1.*6.942540295556983*sin(\t r)});
\draw [shift={(0.6509062280755878,5.739981568676403)}]  plot[domain=5.303501600712838:5.950821915852947,variable=\t]({1.*7.828299241711516*cos(\t r)+0.*7.828299241711516*sin(\t r)},{0.*7.828299241711516*cos(\t r)+1.*7.828299241711516*sin(\t r)});
\draw (1,10) node[anchor=north west] {$1$};
\draw (5,9.3) node[anchor=north west] {$2$};
\draw (8,10) node[anchor=north west] {$3$};
\draw (4.4,11.8) node[anchor=north west] {$4$};
\draw [shift={(5.,15.7)}]  plot[domain=4.227911051347666:5.196866909421713,variable=\t]({1.*6.441273166075166*cos(\t r)+0.*6.441273166075166*sin(\t r)},{0.*6.441273166075166*cos(\t r)+1.*6.441273166075166*sin(\t r)});
\draw [shift={(5.,4.)}]  plot[domain=1.1071487177940904:2.0344439357957027,variable=\t]({1.*6.708203932499369*cos(\t r)+0.*6.708203932499369*sin(\t r)},{0.*6.708203932499369*cos(\t r)+1.*6.708203932499369*sin(\t r)});
\end{tikzpicture}
\end{center}
\caption{Diagramatic representation of the affine cactus $\sigma_{1,2}\sigma_{3,1}\sigma_{2,3}$ dans $AJ_4$}\label{figure1}
\end{figure}
Here, $s_{k,l} \in S_n$ is the permutation that reverses the order in the circular interval $[k,l]_c$ given by 
$$\forall i \in [1,n], ~~ s_{k,l}(i) = \left\{\begin{array}{cl}
k+l-i + n  & \text{if } i \in [k,l]_c \text{ and } i \geq k+l\\
k+l-i  & \text{if } i \in [k,l]_c \text{ and } i < k+l\\
k & \text{otherwise}.
\end{array}\right.$$
The generator $\sigma_{i,j}$ may be represented diagrammatically as a braid with $n$ strands on a cylinder, where the strands $i <_c i+1 <_c \ldots <_c j$ intersect each other in one point and reverse their order after that point, see Figure \ref{figure1}. The relations are illustrated in Figure \ref{figure2}.\\
\begin{figure}[h]
\begin{center}
\begin{tikzpicture}[line cap=round,line join=round,x=.5cm,y=.3cm,rounded corners=3pt,line width=1pt]
\draw (2.,-0.5)-- (2.,-1.4)-- (4.,-2.6)-- (4.,-3.4)-- (2.,-4.6)-- (2.,-5.5);
\draw (3.,-0.5)-- (3.,-1.4)-- (3.,-2.6)-- (3.,-3.4)-- (3.,-4.6)-- (3.,-5.5);
\draw (4.,-0.5)-- (4.,-1.4)-- (2.,-2.6)-- (2.,-3.4)-- (4.,-4.6)-- (4.,-5.5);
\draw (1.,-0.5)-- (1.,-1.4)-- (1.,-2.6)-- (1.,-3.4)-- (1.,-4.6)-- (1.,-5.5);
\draw (5.,-0.5)-- (5.,-1.4)-- (5.,-2.6)-- (5.,-3.4)-- (5.,-4.6)-- (5.,-5.5);
\draw (6.,-0.5)-- (6.,-1.4)-- (6.,-2.6)-- (6.,-3.4)-- (6.,-4.6)-- (6.,-5.5);
\node at (7.5,-3) {$=$ \hspace*{3pt}};	
\end{tikzpicture}
\begin{tikzpicture}[line cap=round,line join=round,x=.5cm,y=.3cm,rounded corners=3pt,line width=1pt]
\draw (1.,-0.5)-- (1.,-5.5);
\draw (2.,-0.5)-- (2.,-5.5);
\draw (3.,-0.5)-- (3.,-5.5);
\draw (4.,-0.5)-- (4.,-5.5);
\draw (5.,-0.5)-- (5.,-5.5);
\draw (6.,-0.5)-- (6.,-5.5);
\node at (7.5,-3) {\hspace*{10pt}};
\end{tikzpicture}
\begin{tikzpicture}[line cap=round,line join=round,x=.5cm,y=.3cm,rounded corners=3pt,line width=1pt]
\draw (1.,-0.5)-- (1.,-1.4)-- (2.,-2.6)-- (2.,-3.4)-- (2.,-4.6)-- (2.,-5.5);
\draw (2.,-0.5)-- (2.,-1.4)-- (1.,-2.6)-- (1.,-3.4)-- (1.,-4.6)-- (1.,-5.5);
\draw (3.,-0.5)-- (3.,-1.4)-- (3.,-2.6)-- (3.,-3.4)-- (6.,-4.6)-- (6.,-5.5);
\draw (4.,-0.5)-- (4.,-1.4)-- (4.,-2.6)-- (4.,-3.4)-- (5.,-4.6)-- (5.,-5.5);
\draw (5.,-0.5)-- (5.,-1.4)-- (5.,-2.6)-- (5.,-3.4)-- (4.,-4.6)-- (4.,-5.5);
\draw (6.,-0.5)-- (6.,-1.4)-- (6.,-2.6)-- (6.,-3.4)-- (3.,-4.6)-- (3.,-5.5);
\node at (7.5,-3) {$=$ \hspace*{3pt}};
\end{tikzpicture}
\begin{tikzpicture}[line cap=round,line join=round,x=.5cm,y=.3cm,rounded corners=3pt,line width=1pt]
\draw (6.,-0.5)-- (6.,-1.4)-- (3.,-2.6)-- (3.,-3.4)-- (3.,-4.6)-- (3.,-5.5);
\draw (5.,-0.5)-- (5.,-1.4)-- (4.,-2.6)-- (4.,-3.4)-- (4.,-4.6)-- (4.,-5.5);
\draw (4.,-0.5)-- (4.,-1.4)-- (5.,-2.6)-- (5.,-3.4)-- (5.,-4.6)-- (5.,-5.5);
\draw (3.,-0.5)-- (3.,-1.4)-- (6.,-2.6)-- (6.,-3.4)-- (6.,-4.6)-- (6.,-5.5);
\draw (1.,-0.5)-- (1.,-1.4)-- (1.,-2.6)-- (1.,-3.4)-- (2.,-4.6)-- (2.,-5.5);
\draw (2.,-0.5)-- (2.,-1.4)-- (2.,-2.6)-- (2.,-3.4)-- (1.,-4.6)-- (1.,-5.5);
\end{tikzpicture}
\hspace*{0.25\textwidth}

\begin{tikzpicture}[line cap=round,line join=round,x=.5cm,y=.3cm,rounded corners=3pt,line width=1pt]
\draw (1.,-0.5)-- (1.,-1.4)-- (6.,-2.6)-- (6.,-3.4)-- (6.,-4.6)-- (6.,-5.5);
\draw (2.,-0.5)-- (2.,-1.4)-- (5.,-2.6)-- (5.,-3.4)-- (5.,-4.6)-- (5.,-5.5);
\draw (3.,-0.5)-- (3.,-1.4)-- (4.,-2.6)-- (4.,-3.4)-- (1.,-4.6)-- (1.,-5.5);
\draw (4.,-0.5)-- (4.,-1.4)-- (3.,-2.6)-- (3.,-3.4)-- (2.,-4.6)-- (2.,-5.5);
\draw (5.,-0.5)-- (5.,-1.4)-- (2.,-2.6)-- (2.,-3.4)-- (3.,-4.6)-- (3.,-5.5);
\draw (6.,-0.5)-- (6.,-1.4)-- (1.,-2.6)-- (1.,-3.4)-- (4.,-4.6)-- (4.,-5.5);
\node at (7.5,-3) {$=$ \hspace*{3pt}};
\end{tikzpicture}
\begin{tikzpicture}[line cap=round,line join=round,x=.5cm,y=.3cm,rounded corners=3pt,line width=1pt]
\draw (6.,-0.5)-- (6.,-1.4)-- (3.,-2.6)-- (3.,-3.4)-- (4.,-4.6)-- (4.,-5.5);
\draw (5.,-0.5)-- (5.,-1.4)-- (4.,-2.6)-- (4.,-3.4)-- (3.,-4.6)-- (3.,-5.5);
\draw (4.,-0.5)-- (4.,-1.4)-- (5.,-2.6)-- (5.,-3.4)-- (2.,-4.6)-- (2.,-5.5);
\draw (3.,-0.5)-- (3.,-1.4)-- (6.,-2.6)-- (6.,-3.4)-- (1.,-4.6)-- (1.,-5.5);
\draw (2.,-0.5)-- (2.,-1.4)-- (2.,-2.6)-- (2.,-3.4)-- (5.,-4.6)-- (5.,-5.5);
\draw (1.,-0.5)-- (1.,-1.4)-- (1.,-2.6)-- (1.,-3.4)-- (6.,-4.6)-- (6.,-5.5);
\end{tikzpicture}
\caption{Examples of relations in affine cactus groups}
\label{figure2}
\end{center}
\end{figure}
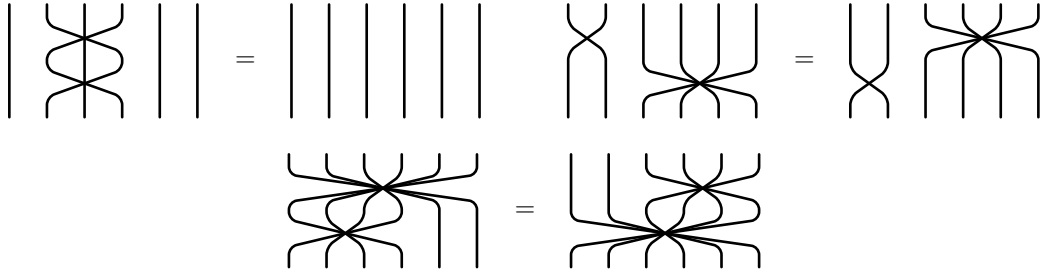
There is a surjective morphism from the affine cactus group $AJ_n$ to the symmetric group $S_n$ given by $$\begin{array}{cccccccccccc}
\pi & : & AJ_n & \to & S_n \\
& & \sigma_{i,j} & \mapsto & s_{i,j}.
\end{array}$$
The kernel of $\pi$ is called the pure affine cactus group and is denoted by $PAJ_n$.\\
Losev \cite{losev_cacti_2019} and Bonnafé \cite{bonnafe_cells_2016} define the notionof generalized cactus group over a Coxeter system $(W,S)$, denoted by $C_W$, as the group generated by $\{\sigma_I, \ \ I \in \mathcal{P}_f^{if}(S)\}$ subject to the following relations:
\begin{eqnarray}\label{Canonical-Presentation}
\sigma_{I}^{2} &=&1 \hspace*{5mm} \textrm{for } I \in \mathcal{P}_f^{if}(S),\\ 
\sigma_{I}\sigma_{J} &=& \sigma_{J}\sigma_{I} \hspace*{5mm} \textrm{if } I \cap J= \emptyset,\\
\sigma_{I}\sigma_{J} &=& \sigma_{J}\sigma_{\omega_J(I)} \hspace*{5mm} \textrm{if } I \subset J.
\end{eqnarray}
Here, $\mathcal{P}_f^{if}(S)$ is the set of subsets $I$ of $S$ for which the parabolic subgroup $W_I$ of $W$ generated by $I$ is finite and irreducible (which corresponds to the fact that the Coxeter diagram of $W_I$ is connected and without cycles \cite{bourbaki_groupes_nodate}) and $\omega_I$ is the longest element in $W_I$ (\cite{bourbaki_groupes_nodate}), acting on the subset $J \in \mathcal{P}_f^{if}(S)$ as follow $$\omega_I(J)=\omega_IJ\omega_I^{-1}.$$
It has been shown by Runze-Yu (\cite{yu_linearity_2023}) that generalized cactus groups are linear, and that the kernel of the morphism $$\begin{array}{ccccccc}
\pi & : & C_W & \to & W \\
& & \sigma_I & \mapsto & \omega_I
\end{array}$$ called the generalized pure cactus group over $W$, and denoted by $PC_W$, embeds in a right-angled Coxeter group.\\
Let $W(\widetilde{A}_n)$ denote the Coxeter group of type $\widetilde{A}_n$. Then we have the following result.
\begin{thmA}
For all $n \geq 2$, we have $$AJ_n \simeq C_{W(\widetilde{A}_n)}.$$
\end{thmA}
To generalise the concept of circular intervals, we introduce the notion of circular sets over $[1,n]$. An ordered set $I=\{i_1 <_c \ldots <_c i_k\}$, where $i_1,\ldots,i_k \in [1,n]$, is called circular if there exists a circular interval $\widetilde{I}$ such that $I$ is contained in $\widetilde{I}$ with preservation of the order. For example, $\{1<_c 3 <_c 7 \}$ is a circular set, but not $\{1 <_c 7 <_c 4 <_c 8\}$. As for circular interval, we say that a circular set $I$ is contained in a circular set $J$, or that $I$ is a c-subset of $J$, denoted by $I \subset_c J$, if the elements of $I$ belong to $J$ with preservation of the order.\\
This allows us to define affine Gauss diagram groups $AD_n$ with generators $\tau_I$, where $I$ is a circular set over $[1,n]$, and whose relations are given by:
\begin{eqnarray}\label{Presentation ADn}
\tau_{I}^{2} &=&1,\\ 
\tau_{I}\tau_{J} &=& \tau_{J}\tau_{I} \hspace*{5mm} \textrm{for } I \cap J= \emptyset \text{ or } I \subset_c J.
\end{eqnarray}
The permutation group $S_n$ acts on the set of circular sets over $[1,n]$ in the following way $$\rho_i(\{i_1 <_c \ldots <_c i_k \})=\left\{\begin{array}{ll}
\{i_1 <_c \ldots <_c i_k \} & \text{if } i,i+1 \in \{i_1,\ldots,i_k\}\\
\{\rho_i(i_1) <_c \ldots <_c \rho_i(i_k) \} & \text{otherwise }
\end{array}\right.$$ where $\rho_i$ is the permutation $(i,i+1)$.\\
This action may seem strange but it allows us to keep track of the notion of circularity. Indeed, it is easy to see that $S_n$ acts by permuting the elements of circular sets then putting them back in circular order. \\
It is obvious that, if $[i,j]_c \subset [k,l]_c$, then $s_{i,j}\cdot \tau_{[k,l]_c} = \tau_{[k,l]_c}$ and $s_{k,l}\cdot\tau_{[i,j]_c}=\tau_{[s_{k,l}(j),s_{k,l}(i)]_c}$. Moreover, if $[i,j]_c \cap [k,l]_c = \emptyset$, then $s_{i,j}\cdot \tau_{[k,l]_c} = \tau_{[k,l]_c}$ and $s_{k,l}\cdot \tau_{[i,j]_c} = \tau_{[i,j]_c}$. Then $S_n$ acts on $AD_n$ and we have the following result.
\begin{thmA}\label{ThB}
For all $n \geq 2$, $AJ_n$ embeds in $AD_n \rtimes S_n$.
\end{thmA}
From this theorem, we can easily deduce that $AJ_n$ is a linear group and that the cactus group $J_n$ embeds in $AJ_n$. Moreover, by considering the restriction of this embedding to pure affine cactus groups, we may see that these groups are residually nilpotent.\\
This embedding  also allows us to study some remarkable subgroups of $AJ_n$. Concretely, for $2 \leq p \leq q \leq n$, we consider the group $AJ_n^{p,q}$ generated by the $\sigma_{i,j}$, where $p \leq [i,j]_c \leq q$, and whose relations are those of $AJ_n$ involving only such generators. In other words, we consider only the generators whose only involve between $p$ and $q$ strands at each crossing. It is easy to check that the groups $AJ_n^{p,p}$ are right-angled Coxeter groups. We have the following result.
\begin{thmA}
For all $2\leq p\leq q \leq n$, the morphism $$\begin{array}{ccccccccccccccc}
\iota & : & AJ_n^{p,q} & \to & AJ_n \\
& & \sigma_{i,j} & \mapsto & \sigma_{i,j}
\end{array}$$ is injective.
\end{thmA}
In particular, the affine cactus groups $AJ_n$ contain the groups $AJ_n^{2,2}$ whose elements may be represented diagrammatically as elements of the Twin groups $Tw_n$ on a cylinder. Moreover, we will deduce the following semi-direct product structure: $$AJ_n \simeq \langle\langle AJ_n^{2,p-1}\rangle\rangle \rtimes AJ_n^{p,n},$$ where $\langle\langle AJ_n^{2,p-1}\rangle\rangle$ is the normal closure of $AJ_n^{2,p-1}$ in $AJ_n$.\\
We will conclude this article with the following theorems, that are consequences of Theorem \ref{ThB} and whose results are similar to those obtained for $J_n$ by Bellingeri-Chemin-Lebed \cite{bellingeri_cactus_2024} for the usual cactus groups.
\begin{thmA}
We have the following results:
\begin{itemize}
    \item The centres of $AJ_n$ and $PAJ_n$ are trivial for $n \geq 2$ and $n \geq 3$ respectively.
    \item The torsion elements of $AJ_n$ are of even order. Moreover, for all $k \geq 1$ there exists an element of order $2^k$ in $AJ_n$ for large enough $n$.
    \item The pure affine cactus group $PAJ_n$ is torsion free for all $n \geq 2$.
    \item The order of the torsion elements in $AJ_n$ is bounded above by $2^{n-1}$.
\end{itemize}
\end{thmA}
\begin{ack}
The author is grateful to John Guaschi and Paolo Bellingeri for their mentoring, helpful insights and careful reading of the paper. The author has received funding from the Normandy region no. 00123353-22E01371.
\end{ack}
\section{Affine cactus groups are generalized cactus groups}
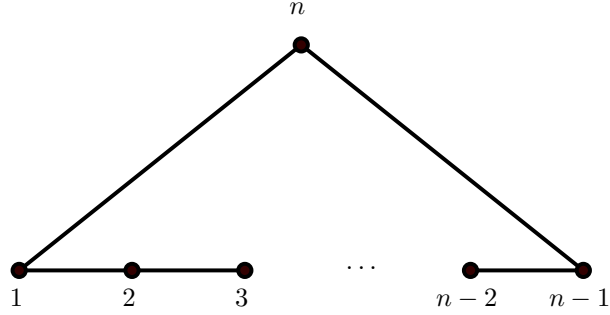
\begin{figure}[h]
\begin{center}
\begin{tikzpicture}[line cap=round,line join=round,x=1.5cm,y=1cm,rounded corners=3pt,line width=1.5pt]
\draw (3.5,3.)-- (1.,0.);
\draw (6.,0.)-- (3.5,3.);
\draw (1.,0.)-- (3.,0.);
\draw (5.,0.)-- (6.,0.);
\draw (3.8,0.2) node[anchor=north west] {$\ldots$};
\draw (0.82,-0.1) node[anchor=north west] {$1$};
\draw (1.82,-0.1) node[anchor=north west] {$2$};
\draw (2.82,-0.1) node[anchor=north west] {$3$};
\draw (4.6,-0.1) node[anchor=north west] {$n-2$};
\draw (5.6,-0.1) node[anchor=north west] {$n-1$};
\draw (3.3,3.7) node[anchor=north west] {$n$};
\begin{scriptsize}
\draw [fill=ttqqqq] (1.,0.) circle (2.5pt);
\draw [fill=ttqqqq] (2.,0.) circle (2.5pt);
\draw [fill=ttqqqq] (3.,0.) circle (2.5pt);
\draw [fill=ttqqqq] (5.,0.) circle (2.5pt);
\draw [fill=ttqqqq] (6.,0.) circle (2.5pt);
\draw [fill=ttqqqq] (3.5,3.) circle (2.5pt);
\end{scriptsize}
\end{tikzpicture}
\end{center}
\caption{Coxeter diagram of type $\widetilde{A}$}\label{figure3}
\end{figure}

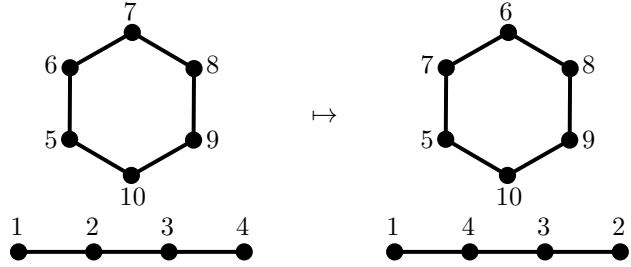
\begin{figure}[h]
\begin{center}
\begin{tikzpicture}[line cap=round,line join=round,x=1cm,y=1cm,rounded corners=3pt,line width=1.5pt]
\draw (1.,1.)-- (2.,1.)-- (3.,1.)-- (4.,1.);
\draw (1.685321236258924,3.4491726483104284)-- (2.5121435458482555,3.9194526757751675)-- (3.332829151319828,3.4385435649874543)-- (3.32669244720207,2.4873544267350027)-- (2.499870137612739,2.0170743992702636)-- (1.6791845321411663,2.497983510057977);
\draw (1.685321236258924,3.4491726483104284)-- (1.6791845321411663,2.497983510057977);
\draw (0.75,1.6) node[anchor=north west] {$1$};
\draw (1.75,1.6) node[anchor=north west] {$2$};
\draw (2.75,1.6) node[anchor=north west] {$3$};
\draw (3.75,1.6) node[anchor=north west] {$4$};
\draw (1.2,2.75) node[anchor=north west] {$5$};
\draw (1.2,3.75) node[anchor=north west] {$6$};
\draw (2.25,4.45) node[anchor=north west] {$7$};
\draw (3.35,3.75) node[anchor=north west] {$8$};
\draw (3.35,2.75) node[anchor=north west] {$9$};
\draw (2.2,2) node[anchor=north west] {$10$};
\draw (4.75,3) node[anchor=north west] {$\Huge\mapsto$};
\draw (6.,1.)-- (7.,1.)-- (8.,1.)-- (9.,1.);
\draw (6.68532,3.44917)-- (7.5121447038257205,3.919447604020764)-- (8.332829407651442,3.438535208041527)-- (8.32668940765144,2.487345208041527)-- (7.49986470382572,2.0170676040207636)-- (6.67918,2.49798);
\draw (6.68532,3.44917)-- (6.67918,2.49798);
\draw (5.75,1.6) node[anchor=north west] {$1$};
\draw (8.75,1.6) node[anchor=north west] {$2$};
\draw (7.75,1.6) node[anchor=north west] {$3$};
\draw (6.75,1.6) node[anchor=north west] {$4$};
\draw (6.2,2.75) node[anchor=north west] {$5$};
\draw (7.25,4.45) node[anchor=north west] {$6$};
\draw (6.2,3.75) node[anchor=north west] {$7$};
\draw (8.35,3.75) node[anchor=north west] {$8$};
\draw (8.35,2.75) node[anchor=north west] {$9$};
\draw (7.2,2) node[anchor=north west] {$10$};
\begin{scriptsize}
\draw [fill=black] (1.,1.) circle (2.5pt);
\draw [fill=black] (2.,1.) circle (2.5pt);
\draw [fill=black] (3.,1.) circle (2.5pt);
\draw [fill=black] (4.,1.) circle (2.5pt);
\draw [fill=black] (1.685321236258924,3.4491726483104284) circle (2.5pt);
\draw [fill=black] (1.6791845321411663,2.497983510057977) circle (2.5pt);
\draw [fill=black] (2.499870137612739,2.0170743992702636) circle (2.5pt);
\draw [fill=black] (3.32669244720207,2.4873544267350027) circle (2.5pt);
\draw [fill=black] (3.332829151319828,3.4385435649874543) circle (2.5pt);
\draw [fill=black] (2.5121435458482555,3.9194526757751675) circle (2.5pt);
\draw [fill=black] (6.,1.) circle (2.5pt);
\draw [fill=black] (7.,1.) circle (2.5pt);
\draw [fill=black] (8.,1.) circle (2.5pt);
\draw [fill=black] (9.,1.) circle (2.5pt);
\draw [fill=black] (6.68532,3.44917) circle (2.5pt);
\draw [fill=black] (6.67918,2.49798) circle (2.5pt);
\draw [fill=black] (7.49986470382572,2.0170676040207636) circle (2.5pt);
\draw [fill=black] (8.32668940765144,2.487345208041527) circle (2.5pt);
\draw [fill=black] (8.332829407651442,3.438535208041527) circle (2.5pt);
\draw [fill=black] (7.5121447038257205,3.919447604020764) circle (2.5pt);
\end{scriptsize}
\end{tikzpicture}
\end{center}
\caption{Action of $\omega_{\{2,3,4,6,7\}}$ on the graph of a Coxeter system $(W,\{1, \ldots, 10\})$}\label{figure5}
\end{figure}

\begin{definition}
Let $n \geq 2$. The affine symmetric group is the Coxeter group $W(\widetilde{A}_n)$ admitting the following presentation:\\
Generators: $\rho_i$, for $1\leq i \leq n$.\\
Relations: \begin{itemize}
    \item $\rho_i^2=1$.
    \item $\rho_i\rho_j=\rho_j\rho_i$ si $|i-j| \geq 2$.
    \item $\rho_i\rho_{i+1}\rho_i=\rho_{i+1}\rho_i\rho_{i+1}$ si $1\leq i \leq n-1$.
    \item $\rho_1\rho_n\rho_1 = \rho_n\rho_1\rho_n$.
\end{itemize}
\end{definition}
Although they are similar, these affine symmetric groups are not the same as those defined by Ilin et al. \cite{ilin_moduli_2023}.\\
It is well known (\cite{bourbaki_groupes_nodate}) that $W(\widetilde{A}_n)$ is the Coxeter group of type $\widetilde{A}_n$ whose Coxeter diagram is depicted in Figure \ref{figure3}.\\
We will show that the affine cactus group is isomorphic to the generalized cactus group over the affine symmetric group. To do this, we start by recalling some definitions and properties of a Coxeter system $(W,S)$, where $S=\{ s_1,\ldots ,s_n\}$ is finite.\\
An interesting fact about Coxeter groups is that they are finite if and only if their Coxeter diagram is of type $A_n$, $B_n$, $D_n$, $E_6$, $E_7$, $E_8$, $F_4$, $G_2$, $H_3$, $H_4$ or $I_p$ (\cite{bourbaki_groupes_nodate}). Another interesting property of Coxeter groups is their irreducibility that is defined below.
\begin{proposition}[\cite{bourbaki_groupes_nodate}]
Let $I \subset S$. Let $W_I$ be the group generated by the elements of $I$ whose relations are those of $W$ that only involve the elements of $I$. Then, the canonical morphism $W_I \to W$ is injective.
\end{proposition}
\begin{definition}
Let $I \subset S$. The group $W_I$ is called irreducible if there do not exist $J,K \subset S$ such that $I=J\cup K$ and $W_I =W_J\times W_K$.
\end{definition}
The irreducibility of $W_I$ may be read off from its Coxeter diagram. Indeed, $W_I$ is irreducible if and only if its Coxeter diagram is connected (\cite{bourbaki_groupes_nodate}).
\begin{definition}
Let $\mathcal{P}_f(S)$ denote the subset $I$ of $S$ such that $W_I$ is finite, and let $\mathcal{P}_{ir,f}(S)$ denote the set of $I \in \mathcal{P}_f(S)$ such that $W_I$ is irreducible.
\end{definition}
\begin{theoreme}[\cite{bourbaki_groupes_nodate}]
Let $I =\{s_{i_1},\ldots, s_{i_k}\} \in \mathcal{P}_{f}(S)$. There exists a unique word of maximal length in $W_I$, denoted by $\omega_I$. The element $\omega_I$ acts on the subsets $J$ of $S$ by conjugacy.
\end{theoreme}
On the level of Coxeter diagrams, $\omega_I$ acts by reversing the order of each point belonging to each connected component of the Coxeter diagram of $W_I$ (see Figure \ref{figure5}).\\
We now define the generalized cactus group over a Coxeter system $(W,S)$.
\begin{definition}[\cite{ilin_moduli_2023},\cite{bonnafe_cells_2016}]
Let $(W,S)$ be a Coxeter system. The generalized cactus group over $W$ is the group $C_W$ with the following presentation :\\
Generators : $\sigma_I$, for $I\in \mathcal{P}_{ir,f}(S)$.\\
Relations : \begin{itemize}
    \item $\sigma_I^2=1$.
    \item $\sigma_I\sigma_J=\sigma_J\sigma_I$, if $W_{IUJ}=W_I\times W_J$.
    \item $\sigma_I\sigma_J=\sigma_J\sigma_{\omega_J(I)}$, if $I \subset J$.
\end{itemize}
\end{definition}
\begin{theoreme}\label{theoreme1}
For all $n \geq 2$, there is an isomorphism $$AJ_n \simeq C_{W(\widetilde{A}_n)}.$$
\end{theoreme}
\begin{proof}
From above, $W=W(\widetilde{A}_n)$ is the Coxeter group of type $\widetilde{A}_n$. We set $S=\{\rho_i\}_{ 1\leq i \leq n}$. Now, if $I=\{\rho_{i_1},\ldots,\rho_{i_k}\} \in \mathcal{P}_{ir,f}(S)$, $(W_I,I)$ is a Coxeter whose diagram is a subdiagram of that of $W$ which is connected and without cycles (since the Coxeter diagrams of type $A_n$, $B_n$, $D_n$, $E_6$, $E_7$, $E_8$, $F_4$, $G_2$, $H_3$, $H_4$ or $I_p$ have none). Thus, the Coxeter diagram of $(W_I,I)$ is of type $A_{|I|}$ as shown in Figure \ref{figure4}, and we deduce that $i_p=i_1+p- 1$ for all $1\leq p \leq k$, which allows us to set $i=i_1$, $i+1=i_2$, ..., $j=i+k-1=i_k$. In particular, $\omega_I$ is the word of $W_I$ which reverses the order of $[i,j]_c$ in the Coxeter diagram of $(W_I,I)$, and it therefore acts on the Coxeter diagram of $(W_I,I)$ in the same way as $s_{i,j}$ acts on $[i,j]_c$.
\begin{figure}[h]
\begin{center}
\begin{tikzpicture}[line cap=round,line join=round,x=1.5cm,y=1cm,rounded corners=3pt,line width=1.5pt]
\draw (1.,0.)-- (2.,0.);
\draw (2.,0.)-- (3.,0.);
\draw (4.,0.)-- (5.,0.);
\draw (3.25,0.2) node[anchor=north west] {$\ldots$};
\draw (0.78,-0.2) node[anchor=north west] {$i_1$};
\draw (1.78,-0.2) node[anchor=north west] {$i_2$};
\draw (2.78,-0.2) node[anchor=north west] {$i_3$};
\draw (3.78,-0.2) node[anchor=north west] {$i_{k-1}$};
\draw (4.78,-0.2) node[anchor=north west] {$i_k$};
\draw (0.78,-0.78) node[anchor=north west] {$i$};
\draw (1.65,-0.78) node[anchor=north west] {$i+1$};
\draw (2.65,-0.78) node[anchor=north west] {$i+2$};
\draw (3.65,-0.78) node[anchor=north west] {$j-1$};
\draw (4.78,-0.78) node[anchor=north west] {$j$};
\begin{scriptsize}
\draw [fill=ttqqqq] (1.,0.) circle (2.5pt);
\draw [fill=ttqqqq] (2.,0.) circle (2.5pt);
\draw [fill=ttqqqq] (3.,0.) circle (2.5pt);
\draw [fill=ttqqqq] (4.,0.) circle (2.5pt);
\draw [fill=ttqqqq] (5.,0.) circle (2.5pt);
\end{scriptsize}
\end{tikzpicture}
\end{center}
\caption{Coxeter diagram of $(W_I,I)$}\label{figure4}
\end{figure}
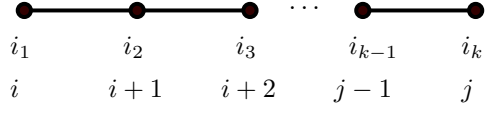
Conversely, for a given $s_{i,j}$, we set $I=\{\rho_i,\ldots,\rho_j\} \in \mathcal{P}_{ir,f}(S)$. The action of $ \omega_I$ on the Coxeter diagram of $(W_I,I)$ corresponds to the permutation $s_{i,j}$ on $[i,j]_c$. Thus, we have the following isomorphism: $$\begin{array}{cccccccc}
S_n & \to & \Omega\\
s_{i,j} & \mapsto & \omega_{\{\rho_i,\ldots,\rho_j\}}
\end{array}$$
where $\Omega$ is the group generated by the $\omega_I$, for $I \in \mathcal{P}_{ir,f}(S)$. Therefore the morphism $$\begin{array}{cccccccc }
AJ_n & \to & J_{AS_n} \\
\sigma_{i,j} & \mapsto & \sigma_{\{\rho_i,\ldots,\rho_j\}}
\end{array}$$ is well defined and it is an isomorphism.
\end{proof}
\section{Linearity and residual nilpotence}
The linearity of $AJ_n$ is proved by Runze and Yu \cite{yu_linearity_2023} and is a consequence of Theorem \ref{theoreme1}. However, we give here an alternative proof from which we also obtain the residual nilpotency of $PAJ_n$.
\begin{proposition}\label{prop1}
Let $n \geq 2$. The map defined by $$\begin{array}{cccccccccccc}
\pi & : & AJ_n & \to & S_n \\
& & \sigma_{i,j} & \mapsto & s_{i,j}
\end{array}$$ is a surjective morphism.
\end{proposition}
\begin{proof}
We just need to check that the relations in $AJ_n$ are sent to relations in $S_n$ via $\pi$.\\
Relations $\sigma_{i,j}^2=1$: Clearly $s_{i,j}^2=1$.\\
Relations $\sigma_{i,j}\sigma_{k,l} = \sigma_{k,l} \sigma_{i,j}$ for $[i,j]_c \cap [k,l]_c = \emptyset$: Since $[i,j]_c \cap [k,l]_c = \emptyset$, the permutations $s_{i,j}$ and $s_{k,l}$ have disjoint supports, and therefore commute.\\
Relations $\sigma_{i,j}\sigma_{k,l}=\sigma_{k,l}\sigma_{s_{k,l}(j),s_{k,l}(i)}$ for $ [i,j]_c\subset_c [k,l]_c$: It suffices to show that $s_{i,j}s_{k,l}(p)=s_{k,l}s_{ s_{k,l}(j),s_{k,l}(i)}(p)$ for all $p \in [1,n]$. Since the permutations of $S_n$ have values in $[1,n]$, it is equivalent to showing that $s_{i,j}s_{k,l}(p)=s_{k,l}s_{ s_{k,l}(j),s_{k,l}(i)}(p) \;(\bmod\; n)$. Now it is easy to verify that these two expressions are equal in $\Z/n\Z$ to: $$\left\{\begin{array}{cl}
p & \text{if } p \notin [k,l]_c \\
k+l-p & \text{if } p \in [k,l]_c \text{ and } p \notin [k+l-j,k+l-i]_c\\
i+j+p -k-l & \text{if } p \in [k,l]_c \text{ and } p \in [k+l-j,k+l-i]_c
\end{array}\right.$$
Thus, $\pi$ is well defined. Furthermore, this morphism is surjective, because the set of transpositions $\{s_{i,i+1}\}_{1\leq i \leq n-1}$ generates $S_n$.
\end{proof}
We will now show that $AJ_n$ embeds in the group $AD_n \rtimes S_n$.\\
\begin{proposition}
Let $n \geq 2$. The morphism $$\begin{array}{cccccccccccccc}
\varphi & : & AJ_n & \to & AD_n \rtimes S_n \\
& & \sigma_{i,j} & \mapsto & (\tau_{[i,j]_c},s_{i,j})
\end{array}$$ is well defined.
\end{proposition}
\begin{proof}
It suffices to check that the relations of $AJ_n$ are sent to relations of $AD_n \rtimes S_n$ via $\varphi$.\\
Relations $\sigma_{i,j}^2=1$: As $s_{i,j}\cdot \tau_{[i,j]_c} = \tau_{[i,j]_c}$, we have: 
$$\left(\tau_{[i,j]_c},s_{i,j}\right)^2 = \left(\tau_{[i,j]_c}s_{i,j}\cdot \tau_{[i,j]_c}, s_{i,j}^2\right) = \left(\tau_{[i,j]_c}^2,s_{i,j}^2\right) = 1.$$
Relations $\sigma_{i,j}\sigma_{k,l} = \sigma_{k,l} \sigma_{i,j}$ for $[i,j]_c \cap [k,l]_c = \emptyset$: Since $[i,j]_c \cap [k,l]_c = \emptyset$, we have $s_{i,j}\cdot \tau_{[k,l]_c} = \tau_{[k ,l]_c}$ and $s_{k,l}\cdot \tau_{[i,j]_c} = \tau_{[i,j]_c}$. Furthermore, $\tau_{[i,j]_c}$ commutes with $\tau_{[k,l]_c}$ and $s_{i,j}$ commutes with $s_{k,l}$. So: $$\begin{array}{lclcl}
\left(\tau_{[i,j]_c},s_{i,j}\right)\left(\tau_{[k,l]_c},s_{k,l}\right) & = & \left(\tau_{[i,j]_c}s_ {i,j}\cdot\tau_{[k,l]_c},s_{i,j}s_{k,l}\right) = \left(\tau_{[i,j]_c}\tau_{[k,l]_c},s_{i,j}s_{k,l}\right)\\
&=& \left(\tau_{[k,l]_c}\tau_{[i,j]_c},s_{k,l}s_{i,j}\right)=\left(\tau_{[k,l]_c}s_{k,l}\cdot\tau_{[i,j]_c},s_{k,l}s_{i,j}\right)\\
&=&\left(\tau_{[k,l]_c},s_{k,l})(\tau_{[i,j]_c},s_{i,j}\right).
\end{array}$$
Relations $\sigma_{i,j}\sigma_{k,l}=\sigma_{k,l}\sigma_{s_{k,l}(j),s_{k,l}(i)}$ for $[i,j]_c\subset_c [k,l]_c$: Like $[i,j]_c\subset_c [k,l]_c$, $s_{i,j}\cdot \tau_{[k,l]_c}=\tau_{[k,l]_c}$ and $s_{k,l} \cdot \tau_{[s_{k,l}(j),s_{k,l}(i)]_c}=s_{k,l}^2 \cdot \tau_{[i,j]_c}=\tau_{[i,j]}$. Furthermore, $\tau_{[i,j]_c}$ and $\tau_{[k,l]_c}$ commute and $s_{i,j}s_{k,l}=s_{k,l}s_{ s_{k,l}(j),s_{k,l}(i)}$ according to Proposition \ref{prop1}. So,
$$\begin{array}{lclcl}
\left(\tau_{[i,j]_c},s_{i,j}\right)\left(\tau_{[k,l]_c},s_{k,l}\right) & = &\left(\tau_{[i,j]_c}s_ {i,j}\cdot\tau_{[k,l]_c},s_{i,j}s_{k,l}\right) = \left(\tau_{[i,j]_c}\tau_{[k,l]_c},s_{i,j}s_{k,l}\right)\\
&=& \left(\tau_{[k,l]_c}\tau_{[i,j]_c},s_{k,l}s_{s_{k,l}(j),s_{k,l}(i) }\right)=\left(\tau_{[k,l]_c}s_{k,l} \cdot \tau_{[s_{k,l}(j),s_{k,l}(i)]_c},s_{k ,l}s_{s_{k,l}(j),s_{k,l}(i)}\right)\\
&=&\left(\tau_{[k,l]_c},s_{k,l}\right)\left(\tau_{[s_{k,l}(j),s_{k,l}(i)]_c},s_{ s_{k,l}(j),s_{k,l}(i)}\right).
\end{array}$$
Thus, the correspondence $\varphi :AJ_n \to AD_n \rtimes S_n$ is a well-defined morphism.
\end{proof}
\begin{lemme}\label{lemma1}
Let $n > 2$, $k\in \N$ and $\sigma =\sigma_{i_1,j_1}\cdots \sigma_{i_k,j_k} \in AJ_n$. We set $\tau=\tau_{I_1}\cdots \tau_{I_k}\in AD_n$, such that $$\varphi(\sigma)=(\tau,s_{i_1,j_1}\cdots s_{i_k, j_k}) \text{ and } I_r=s_{i_1,j_1} \cdots s_{i_{r-1},j_{r-1}}([i_r,j_r]), \text{ for all } r \in [1 ,k].$$
i) If there exists $r \in [2,k]$ such that $\tau_{I_{r-1}}=\tau_{I_r}$, then $\sigma_{i_{r-1}, j_{r-1}}=\sigma_{i_r,j_r}$. \\
In particular, $\sigma = \sigma_{i_1,j_1} \cdots \sigma_{i_{r-2},j_{r-2}}\sigma_{i_{r+1},j_{r+1}} \cdots \sigma_{i_k,j_k}$.\\
ii) If there exists $r \in [2,k]$ such that $\tau_{I_{r-1}}$ and $\tau_{I_r}$ commute, then:
\begin{itemize}
    \item if $I_{r-1} \cap I_r =\emptyset$, then $[i_{r-1},j_{r-1}]_c \cap [i_r,j_r]_c = \emptyset$. \\
    In particular, $\sigma=\sigma_{i_1,j_1} \cdots \sigma_{i_r,j_r}\sigma_{i_{r-1},j_{r-1}} \cdots \sigma_{i_k,j_k}$ .
    \item if $I_{r-1} \subset_c I_r$, then $[i_{r-1},j_{r-1}]_c \subset_c [i_r,j_r]_c$. \\
    In particular, $\sigma=\sigma_{i_1,j_1} \cdots \sigma_{i_r,j_r}\sigma_{s_{i_r,j_r}(j_{r-1}),s_{i_r,j_r}(i_{ r-1})} \cdots \sigma_{i_k,j_k}$.
    \item if $I_r \subset_c I_{r-1}$, then $[i_r,j_r]_c \subset_c [i_{r-1},j_{r-1}]_c$.\\
    In particular, $\sigma=\sigma_{i_1,j_1} \cdots \sigma_{s_{i_{r-1},j_{r-1}}(j_r),s_{i_{r-1},j_{ r-1}}(i_r)}\sigma_{i_{r-1},j_{r-1}} \cdots \sigma_{i_k,j_k}$.
\end{itemize}
\end{lemme}
\begin{proof}
i) Let $r \in [2,k]$ such that $\tau_{I_{r-1}}=\tau_{I_r}$. We then have $$s_{i_1,j_1}\cdots s_{i_{r-2},j_{r-2}}([i_{r-1},j_{r-1}]_c) =_c s_{i_1 ,j_1}\cdots s_{i_{r-1},j_{r-1}}([i_{r},j_{r}]_c).$$
Thus, we have $[i_{r-1},j_{r-1}]_c=_c s_{i_{r-1},j_{r-1}}([i_{r},j_{r}]_c) $.\\ 
Since $s_{i_{r-1},j_{r-1}}$ is an involution, we obtain $s_{i_{r-1},j_{r-1}}([i_{ r-1},j_{r-1}]_c)=_c [i_{r},j_{r}]_c$. \\Now, $s_{i_{r-1},j_{r-1}}([i_{r-1},j_{r-1}]_c)=_c[i_{r-1},j_{r- 1}]_c$, and therefore: $$[i_{r-1},j_{r-1}]_c=_c[i_{r},j_{r}]_c.$$
ii) Let $r \in [2,k]$ such that $\tau_{I_{r-1}}$ and $\tau_{I_r}$ commute. According to the presentation of $AD_n$ (\ref{Presentation ADn}), we then have $I_{r-1}\cap I_r = \emptyset$, or $I_{r-1} \subset_c I_r$ or $I_r \subset_c I_{r-1 }$.\\
Case 1: $I_{r-1}\cap I_r = \emptyset$.\\
Since $s_{i_{r-1},j_{r-1}}([i_{r-1},j_{r-1}]_c)=_c[i_{r-1},j_{r-1} ]_c$, we have:
$$\begin{array}{lcl}
I_{r-1}\cap I_r = \emptyset & \Leftrightarrow & s_{i_1,j_1} \cdots s_{i_{r-2},j_{r-2}}([i_{r-1},j_ {r-1}]_c) \cap s_{i_1,j_1} \cdots s_{i_{r-1},j_{r-1}}([i_{r},j_{r}]_c) = \emptyset \\
& \Leftrightarrow & [i_{r-1},j_{r-1}]_c \cap s_{i_{r-1},j_{r-1}}([i_{r},j_{r}]_c) = \emptyset\\
& \Leftrightarrow & s_{i_{r-1},j_{r-1}}([i_{r-1},j_{r-1}]_c) \cap s_{i_{r-1},j_{ r-1}}([i_{r},j_{r}]_c) = \emptyset\\
& \Leftrightarrow & [i_{r-1},j_{r-1}]_c \cap [i_{r},j_{r}]_c = \emptyset.
\end{array}$$
Case 2: $I_{r-1}\subset_c I_r$.\\
Since $s_{i_{r-1},j_{r-1}}([i_{r-1},j_{r-1}]_c)=[i_{r-1},j_{r-1} ]_c$, we have:
$$\begin{array}{lcl}
I_{r-1}\subset_c I_r & \Leftrightarrow & s_{i_1,j_1} \cdots s_{i_{r-2},j_{r-2}}([i_{r-1},j_{r- 1}]_c) \subset_c s_{i_1,j_1} \cdots s_{i_{r-1},j_{r-1}}([i_{r},j_{r}]_c) \\
& \Leftrightarrow & [i_{r-1},j_{r-1}]_c \subset_c s_{i_{r-1},j_{r-1}}([i_{r},j_{r}]_c) \\
& \Leftrightarrow & s_{i_{r-1},j_{r-1}}([i_{r-1},j_{r-1}]_c) \subset_c s_{i_{r-1},j_{ r-1}}([i_{r},j_{r}]_c)\\
& \Leftrightarrow & [i_{r-1},j_{r-1}]_c \subset_c [i_{r},j_{r}]_c.
\end{array}$$
Case 3: $I_{r}\subset_c I_{r-1}$.\\
Since $s_{i_{r-1},j_{r-1}}([i_{r-1},j_{r-1}]_c)=[i_{r-1},j_{r-1} ]_c$, we have: $$\begin{array}{lcl}
I_{r}\subset_c I_{r-1} & \Leftrightarrow & s_{i_1,j_1} \cdots s_{i_{r-1},j_{r-1}}([i_{r},j_{r }]_c) \subset_c s_{i_1,j_1} \cdots s_{i_{r-2},j_{r-2}}([i_{r-1},j_{r-1}]_c) \\
& \Leftrightarrow & s_{i_{r-1},j_{r-1}}([i_{r},j_{r}]_c) \subset_c [i_{r-1},j_{r-1}]_c \\
& \Leftrightarrow & [i_{r},j_{r}]_c \subset_c s_{i_{r-1},j_{r-1}}([i_{r-1},j_{r-1}]_c )\\
& \Leftrightarrow & [i_{r},j_{r}]_c \subset [i_{r-1},j_{r-1}]_c.
\end{array}$$
This concludes the proof of the lemma.
\end{proof}
\begin{theoreme}\label{theorem2}
The morphism $\varphi: AJ_n \to AD_n \rtimes S_n$ is injective.\\
In particular, $AJ_n$ is a linear group.
\end{theoreme}
\begin{proof}
Let $\sigma =\sigma_{i_1,j_1}\cdots \sigma_{i_k,j_k} \in \ker(\varphi)$ and $\tau=\tau_{I_1}\cdots \tau_{I_k}\in AD_n$ be constructed from $\sigma$ as in Lemma \ref{lemma1}. \\
We then have $\tau = 1$, because $ \sigma \in \ker(\varphi)$, and there therefore exists a sequence of 'movements' of the letters in the word $\tau_{I_1}\cdots \tau_{I_k}$ which allows us to get the empty word. Now, as $AD_n$ is a right-angled Coxeter, the only possible movements are commuting the letters of $\tau$ and their cancellation if the same two letters are next to each other. \\
However, according to Lemma \ref{lemma1}, any cancellation of two letters of $\tau$ corresponds to the cancellation of the two corresponding letters in the expression of $\sigma$, and the commutation of two letters of $\tau$ corresponds to the commutation or 'quasi'-commutation of the two corresponding letters in the expression of $\sigma$.\\
Thus, continuing this process, we deduce that $\sigma=1$ and therefore $\varphi$ is injective.\\
The fact that $AJ_n$ is linear is a consequence of the fact that Coxeter groups are linear and that the semi-direct product of two linear groups is linear (\cite{bourbaki_groupes_nodate}).
\end{proof}
\begin{remark}
In the same way as in Mostovoy's article, we can easily show that the set map $$\begin{array}{ccccc}
\overline{\theta} & : & AJ_n & \to & AD_n \\
& & \sigma_{i,j} & \mapsto & \tau_{[i,j]_c}
\end{array}$$
is injective for $n > 2$.
\end{remark} 
\begin{corollaire}
Pure affine cactus groups $PAJ_n$ are residually nilpotent for $n>2$.
\end{corollaire}
\begin{proof}
According to Theorem \ref{theorem2}, the morphism $\varphi: PAJ_n \to AD_n$ is injective. As $AD_n$ is a right-angled Coxeter group, it is residually nilpotent, and this is also true for all its subgroups, so for $PAJ_n$.
\end{proof}
\begin{remark}
In the case $n=2$, we have $AJ_2=\langle \sigma_{1,2}, \sigma_{2,1} | \sigma_{1,2}^2=\sigma_{2,1}^2=1 \rangle \simeq \Z_2 \ast \Z_2$ and $PAJ_2 = \langle \sigma_{1,2}\sigma_{2,1} \rangle \simeq \Z$. Then $PAJ_2$ is Abelian, and so is nilpotent and residually nilpotent.
\end{remark}.
\section{Direct consequences}
\subsection{Word problem and notable subgroups}
The injectivity of the morphism $\varphi$ gives rise to a certain number of consequences. The first is the following.
\begin{corollaire}\label{corollary2}
The word problem is solvable in $AJ_n$.
\end{corollaire}
\begin{proof}
Let $\sigma=\sigma_{i_1,j_1}\cdots \sigma_{i_k,j_k} \in AJ_n$, and let $\tau = \tau_{I_1}\cdots \tau_{I_k} \in D_n$ be defined as in Lemma \ref{lemma1}. We have $\sigma=1$ if and only if $\tau=1$ by \ref{theorem2}. Now, if $\tau=1$, there exists an algorithm that transforms $\tau$ into the empty word, and which only involves the commutation or the cancelation of two successive letters. Now, according to Lemma \ref{lemma1}, the letters commuting in $\tau$ correspond to the letters in the expression of $\sigma$ which commute or 'quasi'-commute, which provides us with an algorithm to rewrite $\sigma$ as the empty word.
\end{proof}
\begin{definition}
If $2 \leq p \leq q \leq n$, let $AJ_n^{p,q}$ be the group generated by the $\sigma_{i,j}$, such that $p \leq |[i,j ]_c| \leq q$, and subject to the relations of $AJ_n$ involving only these elements.
\end{definition}
The group $AJ_n^{p,q}$ corresponds to cacti which only involve crossings between $p$ and $q$ strands.
\begin{theoreme}\label{theorem3}
The morphism $$\begin{array}{cccccccccccc}
\iota & : & AJ_n^{p,q} & \to & AJ_n \\
& & \sigma_{i,j} & \to & \sigma_{i,j}
\end{array}$$ is injective.
\end{theoreme}
\begin{proof}
Let $AD_n^{p,q}$ denote the group generated by the elements $\tau_I$ belonging to $AD_n$, such that $p \leq |I| \leq q$, subject to the relations of $AD_n$ involving only these elements. \\
As $AD_n^{p,q}$ is a parabolic subgroup of $AD_n$, the inclusion morphism $\iota: D_n^{p,q} \to D_n$ is injective \cite{bourbaki_groupes_nodate}. Moreover, the action of $S_n$ on $AD_n$ restricts to an action on $AD_n^{p,q}$, and as for Theorem \ref{theorem2}, the morphism $$\begin{array}{cccccccc}
\varphi & : & AJ_n^{p,q} & \to & AD_n^{p,q} \rtimes S_n \\
& & \sigma_{i,j} & \mapsto & (\tau_{[i,j]_c},s_{i,j})
\end{array}$$ is injective. \\
Furthermore, we have the following commutative diagram of exact sequences:
$$ \xymatrix{
1 \ar@{->}[r] & AD_n^{p,q} \ar@{^{(}->}[r]\ar@{^{(}->}[d]_\iota & AD_n^{p,q} \rtimes S_n \ar@{->>}[r]\ar@{->}[d]_\alpha & S_n \ar@{->}[r] \ar@{->}[d]_\sim & 1 \\    
    1 \ar@{->}[r] & AD_n \ar@{^{(}->}[r] & AD_n \rtimes S_n \ar@{->>}[r] & S_n \ar@{->}[r]& 1.}$$
According to the five lemma, $\alpha$ is therefore injective, and so we have the following commutative diagram:
$$\xymatrix{
AJ_{n}^{p,q} \ar@{^{(}->}[rd]_{\alpha\circ \varphi} \ar@{->}[r]^\iota & AJ_{n } \ar@{^{(}->}[d]^{\varphi} \\
& AD_{n} \rtimes S_n.
  }$$
Thus, $\iota$ is injective as required
\end{proof}
\begin{remark}
Mostovoy (\cite{mostovoy_round_2023}) defined the full annular twin groups $\overline{B}_n(\mathbb{S}^1)$ by considering the configuration space of $n$ points on a circle, whose fundamental group is the pure annular twin group $\overline{P}_n(\mathbb{S}^1)$, and he showed that $\overline{B}_n(\mathbb{S}^1)$ is generated by the set $\{\sigma_1,\ldots,\sigma_n,\eta\}$ subject to the relations:
\begin{itemize}
    \item $\sigma_i^2 = 1.$
    \item $\sigma_i\sigma_j=\sigma_j\sigma_i$ if $i \not\equiv j \pm 1 \;(\bmod\; n).$
    \item $\sigma_i\eta=\eta\sigma_j$ where $j\equiv i+1 \;(\bmod\; n).$
\end{itemize}
Thus $\overline{B}_n(\mathbb{S}^1) \cong AJ_n^{2,2}\rtimes \Z$, and so it is a subgroup of $AJ_n \rtimes \Z$, where $\Z=\langle \eta \rangle$ acts on $AJ_n$ in the following way: 
$$\eta\sigma_{i,j}\eta^{-1} = \sigma_{p,q} \text{, where } i \equiv p+1 \;(\bmod\; n) \text{ and } j \equiv q+1 \;(\bmod\; n).$$
Moreover, it has been shown that $\overline{P}_n(\mathbb{S}^1)$ is an annular diagram group (\cite{farley_planar_2021}), and that the pure virtual twin group $PVT_n$ embeds in $\overline{P}_n(\mathbb{S}^1)$ (\cite{bellingeri_right-angled_2023}).
\end{remark}
\begin{corollaire}\label{corollary1}
There is a surjective morphism $$\begin{array}{cccclccccccccccc}
\epsilon & : & AJ_n & \to & AJ_n^{p,n}\\
& & \sigma_{i,j} & \mapsto & \left\{\begin{array}{ll}
1 & \text{si } |[i,j]_c| <p\\
\sigma_{i,j} & \text{otherwise.}
\end{array}\right.
\end{array}$$
\end{corollaire}
\begin{proof}
We start by showing that $\epsilon$ is a morphism. To do this, we simply check that the relations of $AJ_n$ are sent to relations of $AJ_n^{p,n}$ via $\epsilon$.\\
Relations $\sigma_{i,j}^2=1$: clearly $\epsilon(\sigma_{i,j})^2=1$.\\
Relations $\sigma_{i,j}\sigma_{k,l}=\sigma_{k,l}\sigma_{i,j}$, for $[i,j]_c\cap[k,l]_c = \emptyset$: without loss of generality, we may assume that $|[i,j]_c| \leq |[k,l]_c|$.
\begin{itemize}
    \item If $|[k,l]_c| < p$, the relation is sent by $\epsilon$ to the trivial relation.
    \item If $|[i,j]_c|<p\leq |[k,l]_c|$, the relation is sent by $\epsilon$ to $\sigma_{k,l}=\sigma_{k, l}$.
    \item If $p\leq |[i,j]_c| \leq|[k,l]_c|$, the relation is sent by $\epsilon$ to $\sigma_{i,j}\sigma_{k,l}=\sigma_{k,l}\sigma_{i, j}$, which is indeed satisfied, because $[i,j]_c\cap[k,l]_c =\emptyset$.
\end{itemize}
Relation $\sigma_{i,j}\sigma_{k,l}=\sigma_{k,l}\sigma_{s_{k,l}(j),s_{k,l}(i)}$, for $[i,j]_c\subset_c [k,l]_c$: we have $|[i,j]|_c \leq |[k,l]|_c$ and $|[s_{k,l}( j),s_{k,l}i]_c|=|s_{k,l}[i,j]_c|=|[i,j]_c|$.
\begin{itemize}
    \item If $|[k,l]_c| < p$, the relation is sent by $\epsilon$ to the trivial relation.
    \item If $|[i,j]_c|<p\leq |[k,l]_c|$, the relation is sent by $\epsilon$ to $\sigma_{k,l}=\sigma_{k, l}$.
    \item If $p\leq |[i,j]_c| \leq|[k,l]_c|$, the relation is sent by $\epsilon$ to $\sigma_{i,j}\sigma_{k,l}=\sigma_{k,l}\sigma_{s_{ k,l}(j),s_{k,l}(i)}$, which is satisfied, because $[i,j]_c\subset [k,l]_c$.
\end{itemize}
Thus, $\epsilon$ defines a morphism. Moreover, it is surjective, because $\epsilon\circ\iota=\Id$, where $\iota: AJ_n^{p,n} \to AJ_n$ is the injective morphism defined in Theorem \ref{theorem3}.
\end{proof}
\begin{corollaire}
There is a split short exact sequence:
\[1 \longrightarrow \langle\langle AJ_n^{2,p-1} \rangle\rangle \overset{\iota}{\longrightarrow} AJ_n \underset{\iota}{\overset{\epsilon}{\myrightleftarrows{ \rule{.7cm}{0cm}}}} AJ_n^{p,n} \longrightarrow 1.\]
In particular, we have: $$AJ_n \simeq \langle\langle AJ_n^{2,p-1} \rangle\rangle \rtimes AJ_n^{p,n}.$$
\end{corollaire}
\begin{proof}
According to Corollary \ref{corollary1}, it suffices to show that $\ker(\epsilon) \subset \langle\langle AJ_n^{2,p-1}\rangle\rangle$, because the inclusion in the other sense is obviously true.\\
Let $\sigma \in \ker(\epsilon)$, and $w$ be a word representing $\sigma$. We may write $w=a_1b_1\cdots a_kb_k$, with $k \in \N$, $a_i \in AJ_n^{p,n}$ and $b_i \in AJ_n^{2,p-1}$ , because $AJ_n$ is generated by the reunion $AJ_n^{p,n}\cup AJ_n^{2,p-1}$, such that all $a_i$ and all $b_i$ are non-empty, except perhaps $a_1$ and $b_k$.\\
We may easily show by induction that: $$\begin{array}{lcl}
w & = & a_1b_1\cdots a_kb_k \\
& = & (a_1b_1a_1^{-1})(a_1a_2b_2a_2^{-1}a_1^{-1})\cdots(a_1 \cdots a_k b_k a_k^{-1} \cdots a_1^{-1})a_1\cdots a_k.
\end{array}$$
Now, since $\sigma \in \ker(\epsilon)$, $a_1\cdots a_k$ represents the trivial element, and therefore $w$ represents the same element as $(a_1b_1a_1^{-1})(a_1a_2b_2a_2^{- 1}a_1^{-1})\ldots(a_1 \ldots a_k b_k a_k^{-1} \ldots a_1^{-1})$, which in turn represents an element of $\langle\langle AJ_n^{2,p- 1}\rangle\rangle$.\\
Moreover, it is clear that $\epsilon \circ \iota = \Id$, and we thus obtain the result.
\end{proof}
We can extend this argument and show that for all $2 \leq r \leq p \leq q \leq m \leq n$, we may viewed $AJ_n^{p,q}$ as a subgroup of $AJ_n^{r,m}$ via the application $\sigma_{i,j}\mapsto\sigma_{i,j}$, and that, for all $2 \leq r < p \leq q \leq n $, we have the following decomposition as a semi-direct product $$AJ_n^{r,q}\simeq \langle\langle AJ_n^{r,p-1}\rangle \rangle \rtimes AJ_n^{p,q}.$$
\subsection{Relations between usual cactus groups and affine cactus groups}
We recall that the usual cactus group of order $n \geq 2$ is the group $J_n$ given by the following presentation:\\
Generators: $\sigma_{i,j}$, for $1 \leq i < j \leq n$.\\
Relations: \begin{itemize}
    \item $\sigma_{i,j}^2=1$.
    \item $\sigma_{i,j}\sigma_{k,l}=\sigma_{k,l}\sigma_{i,j}$ if $[i,j]\cap [k,l]= \emptyset $.
    \item $\sigma_{i,j}\sigma_{k,l}=\sigma_{k,l}\sigma_{k+l-j,k+l-i}$ if $[i,j]\subset [k,l ]$.
\end{itemize}
The group $J_n$ may also be defined as the generalized cactus group over the Coxeter group $S_n$ of type $A_{n-1}$ (\cite{ilin_moduli_2023}).\\
It has been shown (\cite{mostovoy_pure_2019}) that $J_n$ also embeds into $D_n \rtimes S_n$ via the morphism: $$\begin{array}{cccccccc}
\psi & : & J_n & \to & D_n \rtimes S_n \\
& & \sigma_{i,j} & \mapsto & (\tau_{[i,j]},s_{i,j})
\end{array}$$ where $D_n$ is the $n$th Gauss diagram group defined as follows.
\begin{definition}[\cite{mostovoy_pure_2019}]
Let $n \geq 2$. The Gauss diagram group of order $n$ is the group $D_n$ admitting the following presentation:\\
Generators: $\tau_I$, for $I \subset [1,n]$.\\
Relations: \begin{itemize}
    \item $\tau_I^2=1$.
    \item $\tau_I\tau_J=\tau_J\tau_I$ if $I \subset J$ or if $I \cap J =\emptyset$.
\end{itemize}
\end{definition}
The group $S_n$ acts on $D_n$ in the following way: $$\forall s \in S_n,~~ \forall I \subset [1,n],~~ s\cdot \tau_I=\tau_{s(I )}.$$
We thus have the following result.
\begin{theoreme}\label{theorem4}
Let $n \geq 2$. The canonical morphism $\iota: J_n \to AJ_{n}$ is injective.
\end{theoreme}
\begin{proof}
Since $D_n$ is a parabolic subgroup of $AD_n$, we have the following commutative diagram:
$$ \xymatrix{
    1 \ar@{->}[r] & D_n \ar@{^{(}->}[r]\ar@{^{(}->}[d] & D_n \rtimes S_n \ar@{->>}[r]\ar@{->}[d]_\alpha & S_n \ar@{->}[r] \ar@{->}[d]_\sim & 1 \\
    1 \ar@{->}[r] & AD_n \ar@{^{(}->}[r] & AD_n \rtimes S_n \ar@{->>}[r] & S_n \ar@{->}[r]& 1.}$$
According to the five lemma, $\alpha$ is therefore injective, and so is $\alpha \circ \psi : J_n \to AD_n \rtimes S_n$.\\
Then, from Theorem \ref{theorem2}, we have the following commutative diagram:
$$\xymatrix{
J_{n} \ar@{^{(}->}[rd]_{\alpha \circ \psi} \ar@{->}[r]^\iota & AJ_{n} \ar@{^{(}->} [d]^{\varphi} \\
& AD_{n} \rtimes S_{n}.
  }$$
We then have $$\sigma \in \ker(\iota)  \Leftrightarrow  \alpha \circ \psi(\sigma)= \varphi(\iota(\sigma))=1  \Leftrightarrow  \sigma \in \ker(\alpha\circ\psi) \Leftrightarrow \sigma = 1.$$
Thus, the morphism $J_n \to AJ_n$ is injective as required.
\end{proof}
\begin{corollaire}
Let $n \geq 2$. The morphism $PJ_n \to PAJ_n$ is injective.
\end{corollaire}
\begin{proof}
This is an obvious consequence of Theorem \ref{theorem4}.
\end{proof}
\subsection{Torsion in affine cactus groups}
It is known (\cite{bourbaki_groupes_nodate}) that if an element $g$ of a Coxeter group $G$ is central then for any reduced word $w$ representing $g$, each letter of $w$ commutes with any generator of $G$. Then since there is no generator commuting with all of the other generators of $AD_n$, it follows that: $$Z(AD_n)=\{1\}.$$
From this, we deduce the following theorem.
\begin{theoreme}\label{theorem5}
For all $n \geq 2$, the centre of $AJ_n$ is trivial.
\end{theoreme}
\begin{proof}
For the case $n=2$, we have $AJ_2=\langle \sigma_{1,2}, \sigma_{2,1} \ | \ \sigma_{1,2}^2=\sigma_{2,1}^2=1\rangle\simeq \Z/2\Z \ast \Z/2\Z$ and therefore the centre of $AJ_2$ is trivial. We may thus assume that $n \geq 3$.\\
Let $\sigma \in Z(AJ_n) \setminus \{1\}$, and let $w$ be a word of minimal length that represents $\sigma$. Let $\tau \in AD_n$ and $s \in S_n$, be such that $\varphi(\sigma)=(\tau,s)$. We will show that $\tau \in Z(AD_n)=\{1\}$, and therefore $\tau=1$.\\
Let $\sigma_{i,j}$ be a generator of $AJ_n$. Since $\sigma \sigma_{i,j} = \sigma_{i,j}\sigma$, we only have two options according to Corollary \ref{corollary2}.\\
\textbf{Case 1:} The words $\sigma_{i,j}w$ and $w\sigma_{i,j}$ are irreducible. We can then transform $\sigma_{i,j}w$ into $w\sigma_{i,j}$ via commutations and 'quasi'-commutations. In particular, there is a letter $l$ in the expression for $\sigma_{i,j}w$ which can be moved in this way to the end of the word to give the letter $\sigma_{i, j}$. We then have $2$ cases:\\
\textbf{Subcase 1:} The letter $l$ is the initial letter $\sigma_{i,j}$ of $\sigma_{i,j}w$. At the level of the groups of diagrams, according to Lemma \ref{lemma1}, we then have $\tau\tau_{[i,j]_c}=\tau_{[i,j]_c}\tau$. However this means that $\tau$ commutes with all of the generators $\tau_{[i,j]_c}$, and then $\tau=1$.\\
\textbf{Subcase 2:} The letter $l$ appears in the word $w$. Thus, we can transform $w$ into $w'\sigma_{i,j}$ and then $w\sigma_{i,j}=w'$, which contradicts the minimality of $w\sigma_{i,j }$.\\
\textbf{Case 2:} The words $w\sigma_{i,j}$ and $\sigma_{i,j}w$ are simultaneously reducible (since two irreducible representatives of the same element of $AJ_n$ have the same length). As $w$ is minimal, we are in Subcase 2 mentioned above, and therefore there exists a letter $l'$ of $w$ which is sent to the end of $w$ where we transform $w$ into $w' \sigma_{i,j}$. Similarly, there is a letter $l''$ of $w$ which is sent to the beggining of $w$ where we transform $w$ into $\sigma_{i,j}w''$. Once again, two cases arise:\\
\textbf{Subcase 1:} The letters $l'$ and $l''$ occupy the same position in $w$. In particular, all of the letters of $\tau$ commute with $\tau_{[i,j]}$. However this means that $\tau$ commutes with all the generators $\tau_{[i,j]_c}$, and then $\tau=1$.\\
\textbf{Subcase 2:} The letters $l'$ and $l''$ occupy different positions in $w$. We may then transform $w$ into $\sigma_{i,j}u\sigma_{i,j}$. Now, as $\sigma_{i,j}w=w\sigma_{i,j}$, we have $\sigma_{i,j}u=u\sigma_{i,j}$ and therefore $$w= \sigma_{i,j}u\sigma_{i,j}=u\sigma_{i,j}^2=u.$$
Thus $u$ is a word shorter than $w$ representing $\sigma$, which contradicts the minimality of $w$.\\
Thus, for all $1 \leq i \neq j \leq n$, we have $\tau_{i,j}\tau=\tau\tau_{i,j}$ and the only element $\tau$ verifying this equality is $1$. Now, the only element of $AJ_n$ sent to $1$ in $AD_n$ is $1$ by injectivity of $\varphi$ and so the centre of $AJ_n$ is trivial.
\end{proof}
\begin{theoreme}
For $n \geq 3$, the centraliser of $PAJ_n$ in $AJ_n$ is trivial. In particular, the centre of $PAJ_n$ is trivial.
\end{theoreme}
The case $n=2$ is easily handled. Indeed, it is straightforward to show that: $$PAJ_2 =\langle \sigma_{1,2}\sigma_{2,1} \rangle \simeq \Z,$$ and it is therefore an Abelian group.
\begin{proof}
Let $\sigma \in AJ_n$ be a non-trivial element that commutes with all pure affine cacti, and let $w$ be a word of minimum length representing $\sigma$.\\
As $S_n$ is generated by transpositions of the type $(i,i+1)$, with $1\leq i \leq n-1$, any generator $\sigma_{i,j}$ may be transformed into a pure affine cactus defined by: $$\widetilde{\sigma} _{i,j}=\sigma_{i,j}\sigma_{i_1,i_1+1}\sigma_{i_2,i_2+1} \cdots .$$
We then have $\sigma\widetilde{\sigma}_{i,j}=\widetilde{\sigma}_{i,j}\sigma$ by hypothesis and, as in the proof of Theorem \ref{theorem5}, by taking $\tau \in AD_n$ such that $\varphi(\sigma)=(\tau,\pi(\sigma))$, we show that all of the elements $\tau_{[i ,j]}$ commute with all of the letters of $\tau$, and therefore $\tau =1$, which implies that $\sigma = 1$. Thus the centre of $PAJ_n$ in $AJ_n$ is trivial.
\end{proof}
\begin{theoreme}\label{theorem6}
The affine cactus group $AJ_n$ has no odd-order torsion elements.
\end{theoreme}
\begin{proof}
Let $p$ be an odd prime number, let $\sigma \in J_n$ be a non-trivial element of order $p$ such that $\sigma$ has the shortest possible representative $w$. Thus, for every letter $l$ of $w$, there exists a letter $l'$ of $w$, such that $l$ is moved towards $l'$ and cancels with it. We may assume that $l$ moves to the right of the word $w$.\\
\textbf{Case 1:} The letters $l$ and $l'$ occupy different positions in the $q$th and $q'$th copies of $w$ respectively, where $1 \leq q \leq q' \leq p$. Since $w$ is irreducible, we have $q<q'$ We can assume that $l$ is located in the last position and that $l'$ is located in the first position even if it means moving them in $w$. Let $w'$ be the word obtained by sending the first letter of $w$ to the end of this same word, which represents a non-trivial element $\sigma'\in AJ_n$ of order $p$, which is a conjugate of $\sigma$, and where the letters $l$ and $l'$ can still be moved and then canceled. The letter $l$ then belongs to the $q$th copy of $w'$, while $l'$ belongs to the $(q'-1)$th copy of $w'$. By induction, we obtain a non-trivial element $\widetilde{w}$ of order $p$, where the letters $l$ and $l'$ cancel in the same copy $q$ of $w$, which contradicts the minimality of $w$.\\
\textbf{Case 2:} The letters $l$ and $l'$ occupy the same position $i$ in different copies of $w$.\\
We have $\varphi(\sigma)=(\tau,\pi(\sigma))$ and $$1=\varphi(\sigma)^p=(\tau s\cdot\tau s^2\cdot\tau \ldots s^{p-1}\cdot \tau,s^p).$$
Let $u = \tau s\cdot\tau s^2\cdot\tau \cdots s^{p-1}\cdot \tau$, and let $\widetilde{\tau}$ be the letter corresponding to $l$ in $\tau$. The letters of $u$ corresponding to the $p$ copies of the $i$th letter $l$ are $\widetilde{\tau}$, $s\cdot \widetilde{\tau}$, ..., $s ^{p-1}\cdot \tau$. As $p$ is prime, $s$ is of order $p$ or $1$ and in $\{\widetilde{\tau}, s\cdot \widetilde{\tau}, ..., s^{p -1}\cdot \tau\}$, two elements, corresponding to $l$ and $l'$ in $w^p$, can be moved to cancel each other. Thus, the $p$ letters in $u$ corresponding to $l$ are identical and can be moved through the word $u$. Now, as $u$ represents the trivial element of $AD_n$, the letter $\widetilde{\tau}$ appears an even number of times in $u$, and as $p$ is odd, there exists at least one copy of $\widetilde{\tau}$ different from the $p$ copies mentioned above. Thus $\widetilde{\tau}$ appears at least twice in each word $\tau$, $s\cdot\tau$,..., $s^{p-1}\cdot\tau$, where two occurrences can be moved to cancel each other, which contradicts the minimality of $w$.
\end{proof}
\begin{theoreme}
For every $k \in \N^*$, such that $2^k\leq n$, there exists an element of order $2^k$ in $AJ_{n}$.
\end{theoreme}
\begin{proof}
We set $t_1=\sigma_{1,2}$ and, for all $k \in \N^*$, $t_{k+1}=t_{k}\sigma_{1,2^{k+1 }}$. Then $t_1$ is of order $2$ and we have: $$t_{k+1}^2=t_k(\sigma_{1,2^{k+1}}t_k\sigma_{1,2^{k+1}})=t_kt_k',$$
where $t_k'=\sigma_{1,2^{k+1}}t_k\sigma_{1,2^{k+1}}$. We will prove by induction that the order of $t_{k+1}$ is $2^{k+1}$ by assuming that the order of $t_k$ is $2^k$.\\
In $t_k'$ we can move the first letter $\sigma_{1,2^ {k+1}}$ to the end to cancel with $\sigma_{1,2^{k+1}}$ at the end, and therefore $t_k'$ can be written as a word in the $\sigma_{i,j}$, with $i >2^k $. Thus, $t_k$ and $t_k'$ commute and have the same order $2^k$, which implies that $t_{k+1}^2=t_kt_k'$ is of order $2^k$ and therefore $t_{k+1}$ is of order $2^{k+1}.$
\end{proof}
According to Theorem \ref{theorem6}, if $\sigma \in AJ_n$ is a torsion element and $w$ is a word representing $\sigma$, we can rewrite $w$ as a new word $w'$ by (quasi-)commuting its letters in such a way that the number of intersecting strands never increases from top to bottom. Moreover, as the letters of $w'$ (quasi-)commute pairwise, $\sigma$ may necessarily be written in the form $\sigma_{i_1,j_1}\cdots \sigma_{i_k,j_k}$, with $[i_q,j_q]_c \subset [i_p,j_p]_c$ or $[i_q,j_q] \cap [i_p,j_p]=\emptyset$ for $p<q$. In particular, this implies that $|[i_1,j_1]_c|\geq \ldots \geq |[i_k,j_k]_c|$. We will thus say that such an element is decreasing, and $w=\sigma_{i_1,j_1}\cdots \sigma_{i_k,j_k}$ is called a decreasing word representing $\sigma$.\\
Now for two affine cacti $\sigma_1 = \sigma_{i_1,j_1} \cdots \sigma_{i_k,j_k}$ and $\sigma_2 =\sigma_{\widetilde{i_1},\widetilde{j_1}} \cdots \sigma_{\widetilde{i_l},\widetilde{j_l}}$, we say that $\sigma_1$ is disjoint from $\sigma_2$ if for all $1\leq p \leq k$ and $1\leq q \leq l$, $[i_p,j_p]_c \cap [\widetilde{i_q},\widetilde{j_q}]_c=\emptyset$. Obviously, if $\sigma_1$ and $\sigma_2$ are disjoint cactus, then they commute. This allows us to define irreducibly decreasing cacti as decreasing elements $\sigma \in AJ_n$ such that $\sigma$ cannot be written as $\sigma_1\sigma_2$ with $\sigma_1$ and $\sigma_2$ disjoint.
\begin{lemme}\label{lemma2}
Let $\sigma\in AJ_n$. If $\sigma$ is a decreasing cactus, then $\sigma \notin PAJ_n$.
\end{lemme}
\begin{proof}
Let $\sigma=\sigma_{i_1,j_1} \cdots \sigma_{i_k,j_k} \in AJ_n$ be such an element, where $w=\sigma_{i_1,j_1}\cdots \sigma_{i_k,j_k}$ is a decreasing word representing $\sigma$. According to Corollary \ref{corollary2}, we may assume that $\sigma_{i_1,j_1}\cdots \sigma_{i_k,j_k}$ is a reduced word representing $\sigma$. We will argue by induction on $k$ to show that $\pi(\sigma)(i_1) \neq i_1$. If $\sigma$ is not irreducibly decreasing, then there exist $\sigma_1$ and $\sigma_2$ disjoint such that $\sigma=\sigma_1\sigma_2$. However, since $\sigma_1$ and $\sigma_2$ are disjoint, then $\pi(\sigma_1)$ and $\pi(\sigma_2)$ have disjoint support in $S_n$ and only one of them acts on $i_1$. Thus we can assume that $\sigma$ is irreducibly decreasing.\\
As all generators are involutions, without loss of generality we may assume that there does not exist $l \in [1,k-1]$, such that $[i_l,j_l]_c=[i_{l+ 1},j_{l+1}]_c$, which implies that $|[i_1,j_1]_c| > |[i_2,j_2]_c| \geq \ldots \geq |[i_k,j_k]_c|$. Using this and the fact that $\sigma$ is a decreasing cactus, $i_1$ or $j_1$ does not belong to the circular intervals $[i_p,j_p]_c$ for $2\leq p \leq k$. Without loss of generality, we may assume that it is $i_1$.\\
If $k=1$, we have $\pi(\sigma)(i_1) = j_1 \neq i_1$.\\
Let us assume that the conclusion is true for words of length $k$ and prove it for those of length $k+1$. \\
Let $\sigma=\sigma_{i_1,j_1}\cdots \sigma_{i_{k+1},j_{k+1}} \in AJ_n$. We set $\widetilde{\sigma}=\sigma_{i_1,j_1}\cdots \sigma_{i_k,j_k}$. Then:
$$\pi(\sigma)(i_1)=s_{i_{k+1 },j_{k+1}}\pi(\widetilde{\sigma})(i_1).$$
By the induction hypothesis, we have $\pi(\widetilde{\sigma})(i_1)\neq i_1$ and so we have two cases.\\
\textbf{Case 1:} If $\pi(\widetilde{\sigma})(i_1) \in [i_{k+1},j_{k+1}]_c$, then $s_{i_{k+ 1},j_{k+1}}\pi(\widetilde{\sigma})(i_1) \in [i_{k+1},j_{k+1}]_c$. Now, $i_1 \notin [i_{k+1},j_{k+1}]_c$ and therefore $\pi(\sigma)(i_1) \neq i_1$.\\
\textbf{Case 2:} If $\pi(\widetilde{\sigma})(i_1) \notin [i_{k+1},j_{k+1}]_c$, then $s_{i_{k+ 1},j_{k+1}}\pi(\widetilde{\sigma})(i_1) =\pi(\widetilde{\sigma})(i_1)\neq i_1$.\\
Thus, $\sigma \notin PAJ_n$ and so decreasing cacti do not belong to $PAJ_n$.
\end{proof}
\begin{theoreme}
For $n \geq 2$, $PAJ_n$ is torsion-free.
\end{theoreme}
\begin{proof}
According to Theorem \ref{theorem6}, it suffices to show that $PAJ_n$ does not admit torsion elements of order $2$. Suppose on the contrary that there exists such a non-trivial element $\sigma \in PAJ_n$, and let $w$ be a minimal length word representing $\sigma$. As in the proof of Theorem \ref{theorem6}, there exists $\tau \in AD_n$, such that $\varphi(\sigma)=(\tau,\pi(\sigma))$ and such that $\tau $ is a product of generators that commute pairwise. In particular, according to Theorem \ref{theorem6}, we can rewrite $w$ into a new word $w'$ by (quasi-)commuting its letters in such a way that the number of intersected strands never increases from top to bottom. Moreover, as the letters of $w'$ (quasi-)commute pairwise, $\sigma$ is a decreasing cactus. Now, according to Lemma \ref{lemma2} such an element cannot be pure. Thus, there is no torsion element in $PAJ_n$.
\end{proof}
We will now show that the order of torsion elements in $AJ_n$ is bounded above.
\begin{theoreme}\label{theorem7}
Let $\sigma \in AJ_n$ a torsion element. The order of $\sigma$ is bounded above by $2^{n-1}$.
\end{theoreme}
\begin{proof}
Let $\sigma \in AJ_n$ be a torsion element of order $k \geq 2$. According to Theorem \ref{theorem6} this element is a decreasing cactus. If $\sigma$ is not irreducibly decreasing, then there exist two decreasing disjoint cacti $\sigma_1$ and $\sigma_2$ such that $\sigma=\sigma_1\sigma_2$. So: $$1=\sigma^k=\sigma_1^k\sigma_2^k$$ and since $\sigma_1$ and $\sigma_2$ are disjoint we have $$1=\sigma_1^k=\sigma_2^k.$$ Thus we can assume that $\sigma$ is irreducibly decreasing.\\
Moreover, if $w=\sigma_{i_1,j_1}\ldots\sigma_{i_k,j_k}$ is a decreasing word representing $\sigma$, as in Lemma \ref{lemma2}, we may assume that $w$ is reduced, $|[i_1,j_1]_c| > |[i_2,j_2]_c| \geq \ldots \geq |[i_k,j_k]_c|$ and $i_1$ does not belong to the circular intervals $[i_p,j_p]_c$ for $2\leq p \leq k$.\\
For $1\leq p \leq n$, let $AJ_n\setminus\{p\}$ denote the group generated by $\{\sigma_{i,j}, \ \ p \notin [i,j]_c\}$ with the relations of $AJ_n$ involving only those elements. From the bijection: $$s:k \mapsto \left\{\begin{array}{ll}
k-p + n & \text{if } 1 \leq k \leq p-1 \\
k-p & \text{if } p + 1 \leq k \leq n
\end{array}\right.$$
we have the isomorphism: $$\begin{array}{cccccccccccccccc}
\rho & : & AJ_n \setminus\{p\} & \to & J_{n-1}\\
& & \sigma_{i,j} & \mapsto & \sigma_{s(i),s(j)}.
\end{array}$$
So, setting $\widetilde{\sigma}=\sigma_{i_2,j_2}\ldots\sigma_{i_k,j_k}$, from Theorem \ref{theorem4} we have: $$\widetilde{\sigma}\in AJ_n\setminus \{i_1\} \hookrightarrow AJ_{n-1}.$$
Now, from Theorem \ref{theorem6}, there exists $s\geq 1$ such that $k=2^s$. Then, $$1=\sigma^{2^s}=(\sigma_{i_1,j_1}\widetilde{\sigma})^{2^s}=(\sigma_{i_1,j_1}\widetilde{\sigma}\sigma_{i_1,j_1}\widetilde{\sigma})^{2^{s-1}}.$$
Since $\sigma$ is a decreasing cactus, we have $\sigma_{i_1,j_1}\widetilde{\sigma}\sigma_{i_1,j_1}=\hat{\sigma} \in AJ_n\setminus\{p\}$ for a certain $1\leq p \leq n$, and thus: $$1=(\hat{\sigma}\widetilde{\sigma})^{2^{s-1}}.$$
Since $\widetilde{\sigma}$ is decreasing, the same is true for $\hat{\sigma}$. Moreover, since $\hat{\sigma}\widetilde{\sigma}$ is a torsion element, it is decreasing. We then have two cases :\\
\textbf{Case 1:} The cacti $\widetilde{\sigma}$ and $\hat{\sigma}$ are disjoint and commute. So we have: $$1=\widetilde{\sigma}^{2^{s-1}}\hat{\sigma}^{2^{s-1}}$$
and as before, we have: $$1=\widetilde{\sigma}^{2^{s-1}}=\hat{\sigma}^{2^{s-1}}.$$
But since $\widetilde{\sigma}$ and $\hat{\sigma}$ may be viewed as elements of $AJ_{n-1}$, we have $$o(AJ_n) \leq 2o(AJ_{n-1})$$ where $o(AJ_n)$ denotes the maximal order of torsion elements in $AJ_n$.\\
\textbf{Case 2:} The cacti $\widetilde{\sigma}$ and $\hat{\sigma}$ are not disjoint. Since $\overline{\sigma}=\hat{\sigma}\widetilde{\sigma}$ is of finite order, it is irreducibly decreasing, and using the fact that $|[i_1,j_1]_c| > |[i_2,j_2]_c| \geq \ldots \geq |[i_k,j_k]_c|$, we have $n-1 \geq |[i_2,j_2]_c| \geq \ldots \geq |[i_k,j_k]_c|$ and $\overline{\sigma}\in AJ_n\setminus\{p\}$, for $1 \leq p \leq n$, then $\overline{\sigma}$ can be viewed as an element of $AJ_{n-1}$. Now, as for Case 1, we have: $$o(AJ_n) \leq 2o(AJ_{n-1}).$$
The fact that $AJ_2 = \langle \sigma_{1,2},\sigma_{2,1} | \sigma_{1,2}^2=\sigma_{2,1}^2=1\rangle \cong \Z_{2}\ast \Z_{2}$ gives us $o(AJ_2)=2$, which concludes the proof.
\end{proof}
\begin{remark}
Since $J_n$ embeds into $AJ_n$ by Theorem \ref{theorem4}, Theorem \ref{theorem7} applies to $J_n$, and therefore $o(J_n)\leq 2^{n-1}$.
\end{remark}

\end{document}